# Study of a measure of efficiency as a tool for applying the principle of least effort to the derivation of the Zipf and the Pareto laws

A. El Kaabouchi[3], F.X. Machu[1], J. Cocks[1], R. Wang[1], Y.Y. Zhu[4] and Q.A. Wang[1,2] *

[1]Laboratoire SCIQ, ESIEA, 9 Rue Vésale, 75005 Paris, France
[2]IMMM, CNRS UMR 6283, Le Mans Université, Ave. O. Messiaen, Le Mans, France
[3]ESTACA, Parc universitaire Laval-Changé, Rue Georges Charpak, Laval, France
[4]College of Mathematics and Computer Science, Wuhan Textile University, Wuhan, China

**Abstract**

The principle of least effort is believed to be a universal rule for living systems. Its application to the derivation of the power law probability distributions of living systems has long been challenging. Recently, a measure of efficiency was proposed as a tool of deriving Zipf's and Pareto's laws directly from the principle of least effort. The present work is a further investigation of this efficiency measure from a mathematical point of view. The aim is to get further insight into its properties and usefulness as a metric of performance. We address some key mathematical properties of this efficiency such as its sign, uniqueness and robustness. We also look at the relationship between this measure and other properties of the system of interest such as inequality and uncertainty, by introducing a new method for calculating non-negative continuous entropy.

*Corresponding author: alexandre.wang@esiea.fr

# 1) Introduction

The idea of achieving more by doing less is common. It is seen in most if not all living systems. Ferrero referred to this rule as the principle of least effort (PLE), and used it to interpret the mental inertia on the basis of the observation that human beings do not like work [1]. Later on Zipf applied the same rule in a quantitative study of word frequency distribution, and wrote: "*The power laws in linguistics and in other human systems reflect an economical rule: everything carried out by human being and other biological entities must be done with* least effort (at least statistically)" [2][3].

This statement of Zipf extends Ferrero's PLE for human being to other biological entities. According to a common understanding [4], the biological entities here can include human beings, animals, insects and even smaller biological agents. There has been considerable efforts to verify the validity of PLE for some of these living entities or agents (see [4] for a review and the references therein). In what follows, we will use the term "living systems" to represent these living agents as well as the systems used or driven by these living agents (linguistic, social, economic, educational, communicational systems etc.). The words "everything carried out" will be referred to as "achievement" of the living systems trying to achieve something with effort.

The most interesting aspect of Zipf's idea is the possibility of deriving the power laws from PLE. This perspective has inspired a lot of efforts to get power laws by minimizing effort [5][6]. The derivation of power laws of living systems directly from PLE has long been challenging. The main obstacle is the mathematical definition of the key quantities such as achievement and effort, both being process-dependent and difficult to define in a general way (see detailed explanation below).

In a recent work [14], we proposed the idea to derive the Zipf and the Pareto power laws (see section 3) from PLE by using the principle of maximum efficiency (PME). The efficiency was defined as the ratio of achievement to effort, similar to the efficiency of thermal engine defined by the ratio of the output work (=achievement) to the input heat (=effort). The general formula of the efficiency was found with the help of the nonadditive property of the efficiency of thermal engines in thermodynamics. A brief review of that approach is given below. In [14], it is shown that the value of the efficiency can easily be calculated for a system if the distribution laws of achievement of the system are known. This opens the way for this efficiency functional to be used as a measure of the performance of the system. One just needs to establish the distribution laws of achievement by either empirical study or numerical simulation and to calculate the efficiency using those laws. This possibility has motivated the present work which is aimed at getting further insight into this metric of performance by studying its mathematical properties.

In the second section, we give a brief description of the context and the motivations of this work as well as our previous work [14], in order to explain why it is necessary to implement PLE by using PME, and to define the efficiency without using achievement and effort. The third section is a summary of our previous work presented in [14]. The reader will see how the efficiency functional comes out of a general nonadditivity of the efficiency of thermal engines, and how it can be maximized to give the Zipf and the Pareto power laws.

In the fourth section, we present the proof for the uniqueness of the efficiency functional with the given efficiency nonadditivity. This uniqueness is the first guarantee of the reliability and the applicability of the functional as a meaningful tool to derive power laws. In the fifth section we study the robustness (or Lesche stability) of the efficiency functional against small variations of probability, which is another guarantee of the applicability of the efficiency functional as a reliable metric of performance.

The sixth section presents a study of the sign and the numerical range of the efficiency, in order to explore the possibility and the usefulness of the efficiency which can be either negative or greater than one. In the seventh section, we study the relationship between the efficiency and the inequality measured by the Gini coefficient $G$ which is an important property of many systems [23]. The eighth section is a study of the relationship between efficiency and informational entropy. We use a generalized definition of informational entropy called varentropy [24][26][28], which provides a more reliable and optimal measure of probabilistic uncertainty avoiding negative informational entropy[28,30]. The last section is the conclusion.

## 2) From principle of least effort to principle of maximum efficiency

The rule of least effort for the living is very appealing. Ferrero supported this rule by relating the activity of the human brain to the laws of inertia in physics and chemistry [1]. The idea of least effort is so intuitive and obvious for living systems that you can hardly say the contrary. In our opinion, achieving more by doing less should hold a crucial place in the survivability of any species, especially when its possession of the resources to use (force, money, time etc.) and the quantity of achievable things it desires are both limited. We can easily imagine the fate of a species in the ruthless competition of the living world if it always tries to spend as much as possible its force and energy in order to achieve as few as possible its food, or to protect as little as possible itself against dangers.

PLE is especially interesting in the perspective of being applied to derive empirical laws in the same way as many variational principles in physics (stationary action, least time, maximum entropy etc.). For this purpose, a functional of effort is needed. However, as mentioned above, the nature of the effort, just as achievement, depends on the system and process involved. An effort can be the expenditure of any resource, such as (mechanical or biological) energy, time, information, money, and even very abstract things such as spiritual, mental or mindful activities. Similarly, an achievement can be anything a living system desires to carry out or to obtain with effort, such as food, income, wealth, information, education, garment, pleasure, honor, city population, firm size, frequency of words, and so forth (further explanation is given below). Obviously, it is very difficult to use a single general mathematical formula for quantifying these quantities, some of which are even non quantifiable, such as pleasure, honor, mental effort etc. this difficulty of defining effort and achievement is, in our opinion, the main obstacle to the application of PLE using variational calculus minimizing effort.

It was this difficulty which led us to the idea of implementing PLE using an equivalent variational calculus referred to as principle of maximum efficiency (PME) to derive the Zipf and the Pareto laws [14]. PME can be regarded as a consequence of PLE for the following reasons. PLE implies either minimization of effort for given achievement (buying a given amount of goods with minimum of money for instance)[1] or maximization of achievement for a given effort (doing maximum of work for a given amount of time, for example). Both these two sides of PLE imply PME because the efficiency is defined as the ratio of achievement to effort. From these considerations, we can say that PLE and PME are two sides of the same mechanism. One of the advantages of PME is the possibility to define a functional to quantify the efficiency without using well-defined achievement and effort. The reader will see in the next section how this is possible [14].

When discussing the implementation of PLE for deriving power laws, it is worth mentioning the work of Mandelbrot [5] and Cancho *et al* [6] who defined functionals of effort from an informational point of view in order to derive the Zipf law for languages from PLE. However, the expression of effort (cost) of [5] was proportional to the log (ln) of rank, and seemed an ad hoc proposition aimed at obtaining the Zipf power law from the variational calculus. Cancho et al [6] proposed a functional of effort proportional to the information communicated between speaker and listener. During communication however, the information transferred should be a product of the effort instead of the effort itself. The effort in this case should be process-dependent[2]. In addition, the approaches of [5] and [6] are limited to language and cannot be applied to other living systems.

---

[1] The principle of good enough is a very good example of PME and PLE : minimizing the effort for a given achievement that is good enough for the *necessary* needs. See [35] and references therein.

It is also worth mentioning a series of constructive efforts to derive power laws by variational calculus within the frameworks of several generalized statistics including the non-extensive statistics based on Thallis and Rényi entropies [7]-[10], κ-statistics [11], superstatistics and others [10][12]. The main idea here is to extend the Boltzmann-Gibbs statistical mechanics which has the conventional logarithmic entropy and exponential probability, to more general formalisms allowing the derivation of power law probability distributions by maximizing entropy [7]-[11] or by considering the effect of thermal fluctuation [12]. It is worth noticing that this development of statistics theory is mainly guided by the principle of maximum informational entropy [13], and not directly associated with PLE for living systems.

### 3) The functional of efficiency and the Zipf and Pareto laws

In order to apply PME, we must define a functional quantifying efficiency. The idea of our previous work is roughly the following: a living agent, or a system composed of living agents, can be regarded as an engine spending heat $Q$ (effort/input) and providing mechanical work $W$ (achievement/output) with an efficiency defined by $\eta = \frac{W}{Q}$ [14]. As mentioned in the above section, the achievement $W$ can be anything the living agents desire to obtain or fulfill, such as food, income, wealth, city population, firm size, frequency of events and words, information and so forth. Similarly, the effort $Q$ can be anything a living agent consumes in order to achieve what he desires, such as physical or mental effort, energy, time, materials, money, information, etc. Now, as discussed above, both $W$ (achievement) and $Q$ (effort) are difficult to be defined quantitatively. We have to define the efficiency $\eta = \frac{W}{Q}$ without using predefined $W$ and $Q$. It sounds strange. But a solution does exist [14] and is summarized below.

---

[2] If you are buying information, your effort is money. If you are searching for information in documents, your effort can be the time spent, and so on. Obviously, the quantity of information you obtain is not necessarily proportional to your effort.

It is well known that the efficiency in thermodynamics is not an additive quantity. Suppose we have two Carnot engines A and B. The engine A absorbs an energy $Q_1$ from the hot heat bath 1, delivers a quantity of work $W_A$, and loses an energy $Q_2$ to the intermediate heat bath 2 ($W_A = Q_1 - Q_2$), leading to the efficiency $\eta(A) = \frac{W_A}{Q_1} = 1 - \frac{Q_2}{Q_1}$. The engine B absorbs an energy $Q_2$ from the heat bath 2, delivers a work $W_B$, and loses an energy $Q_3$ to a cold heat bath 3 ($W_B = Q_2 - Q_3$), giving the efficiency $\eta(B) = \frac{W_B}{Q_2} = 1 - \frac{Q_3}{Q_2}$. The overall efficiency $\eta(C)$ of the ensemble engine $C$ composed of the two engines $A$ and $B$ ($C=A+B$) is defined by

$$\eta(C) = \frac{W_A + W_B}{Q_1} = 1 - \frac{Q_3}{Q_1}$$

since the composite engine C absorbs $Q_1$ from the heat bath 1, loses $Q_3$ to the heat bath 3 and delivers a work $W_C = W_A + W_B = Q_1 - Q_3$. It is straightforward to calculate

$$\eta(C) = \eta(A) + \eta(B) - \eta(A)\eta(B). \tag{1}$$

Lord Kelvin rewrote Eq.(1) as $1 - \eta(C) = (1 - \eta(A))(1 - \eta(B))$. This form of the efficiency nonadditivity helped Kelvin to discover the absolute temperature, which in turn led to the discovery of the second law of thermodynamics about 170 years ago[3].

Notice that the Carnot engine is an ideal machine functioning in reversible process without loss of energy. Real engines always lose energy by friction, vibration or thermal radiation in an irreversible process. A real engine cannot transform all the heat absorbed into work according to the ideal relationship $Q_1 - Q_2 = W_1$. We can introduce in this case a loss coefficient *a* so that Eq.(1) becomes

$$\eta(C) = \eta(A) + \eta(B) + a\eta(A)\eta(B). \tag{2}$$

or $1 + a\eta(C) = (1 + a\eta(A))(1 + a\eta(B))$. The reversible case of the Carnot engine Eq.(1) is recovered when $a = -1$.

Now if a large number of real engines (living agents) functioning randomly altogether form a single big system, and all the agents in the ensemble are making efforts to achieve something

---

represented by a random variable $X$ having $w$ discrete values $x_i$ with $i = 1,2,...w$. The more of $X$ they achieve, the larger is each $x_i$. This quantity $X$ can be any achievement such as food, income, wealth, city population, frequency of events etc. as mentioned in the section 2. At equilibrium (or stationary) states of the whole systems, the agents are distributed over the whole range of $X$ ($x_i$ for all $i$) with $n_i$ agents achieving the value $x_i$. We have $\sum_{n=1}^{w} n_i = N$. The probability $p_i$ of finding an agents achieving $x_i$ is $p_i = \frac{n_i}{N}$, with the normalization $\sum_{i=1}^{w} p_i = 1$.

Due to the statistical nature of this model and the large number of agents distributed over all the values of $X$, it is reasonable to suppose that the total efficiency $\eta_i$ of the agents obtaining the value $x_i$ depends on the number $n_i$ where $\eta_i = f(n_i)$ or $\eta_i = f(p_i)$. The average efficiency $\eta$ of the whole system reads $\eta = \sum_{i=1}^{w} p_i \eta_i$.

Now let us separate the whole ensemble of agents into two independent subsystems $A$ and $B$, having the efficiency $\eta_k(A)$ and $\eta_j(B)$, respectively. The probability distribution of the agents in $A$ is $p_k(A)$ and that in $B$ is $p_j(B)$. The probability distribution $p_{kj}(C)$ of the whole ensemble can be written as

$$p_{kj}(C) = p_k(A)p_j(B) \qquad (3)$$

We choose Eq.(2) to describe the efficiency nonadditivity. This implies a total efficiency given by

$$\eta_{kj}(C) = \eta_k(A) + \eta_j(B) + a\eta_k(A)\eta_j(B) \qquad (4)$$

or $(1 + a\eta_{kj}(C)) = [1 + a\eta_k(A)][1 + a\eta_j(B)]$.

It can be proved (proof given in the section below) that Eq.(3) and Eq.(4) uniquely leads to $(1 + a\eta_i) = p_i^b$ or $\eta_i = \frac{p_i^b - 1}{a}$. The parameter $b$ is related to $a$ for the following reasons. First, due to the fact that the efficiency $\eta_i$ is positive and that $p_i$ is smaller than unity, $b$ should have opposite sign to $a$. Secondly, from Eq.(4), as $a \to 0$, the efficiency tends to the additive limit $\eta_i \to \eta_k(A) + \eta_j(B)$. Taking into account Eq.(3), we expect the limit $\eta_i \to \frac{b\ln p_i}{a} \to -\ln p_i$, implying $b \to -a$ for $a \to 0$. This implies a relationship $b = f(a)$ with any odd function $f(a)$ obeying $\lim_{a \to 0} f(a) \to -a$. However, as can be seen in the variational calculus using $\frac{p_i^b - 1}{a}$, only $b$ will be associated with the power law distributions, $a$ being absorbed in the normalization constant. Therefore, different relationships $b = f(a)$ do not affect the power law derivation.

With no loss of generality, we choose the simplest relation $b = -a$. The average efficiency of the whole ensemble of $N$ agents reads

$$\eta = \sum_{k=1}^{w} p_i \eta_i = \frac{\sum p_i^{1-a} - 1}{a} \quad (5)$$

or $\eta = \frac{1 - \sum p_i^{1+b}}{b}$. The concave property of this formula is shown in Figure 1 for $w = 2$ and several values of $a$ in the interval $0 \leq a \leq 1$. $\eta$ is a monotonically increasing function of $a$ and concave in this domain with its extremum at $p = \frac{1}{2}$. $\eta$ becomes a constant independent of the probability distribution as $a \to 1$. But this case does not happen because $a$ will be limited to $0 < a < 0.5$ for mathematical reasons [14].

In our previous work [14], the efficiency Eq.(5) was maximized as an application of PLE. Since efficiency is defined as the ratio of achievement to effort, its maximum implies either maximization of achievement for given effort or minimization of effort for given achievement. So maximum efficiency must be associated with maximum achievement represented by the average $\bar{X} = \sum p_i x_i$. In other words, the maximum of $\eta$ and the maximum of $\bar{X}$ are mutually conditioned. So the functional $L = \eta + c\bar{X}$ should be maximized in the calculus of variation $\delta(\eta + c\bar{X}) = 0$, where $c$ is a positive multiplier. This calculus of PLE straightforwardly leads to the Pareto law [15] $P(X > x) = \left(\frac{x_m}{x}\right)^\beta$ where $X$ represents the achievement (income for example), $P(X > x)$ the probability of finding a person with an income larger than the value $x$, $x_m$ the smallest income and $\beta = -(\frac{1}{b} + 1) = \frac{1}{a} - 1$ the constant characterizing the Pareto law.

In the same way [14], we have derived the Zipf law $x_r = \frac{x_1}{r^\alpha}$ where $r$ is the rank of the achievement $x_r$ (frequency of words in a text for instance), $\alpha = a\gamma$ and $\gamma$ is the constant characterizing the relationship between the derivative of the Pareto law $p'(X > x)$ and the rank $r$ of the value $x : p'(X > x) \propto r^\gamma$. The reader is referred to [14] for the details of the variational calculus.

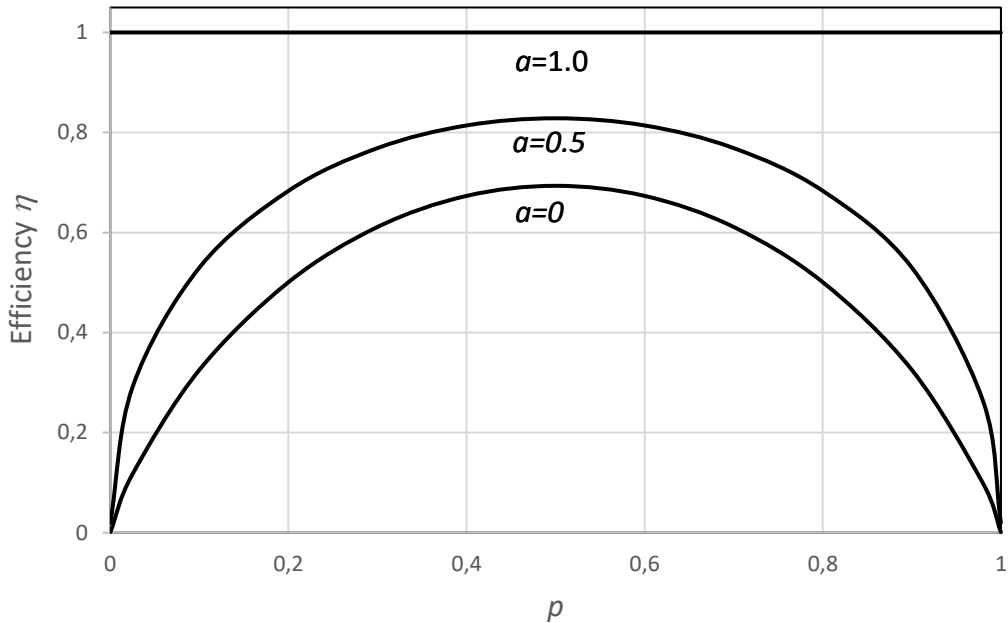

Figure 1. Variation of the efficiency as a function of the probability distribution over two states $p_1 = p, p_2 = 1 - p$ for several values of $a$ in the interval $0 \leq a \leq 1$. $\eta$ is a monotonically increasing function of $a$ in this domain. It is concave with an extremum at $p = \frac{1}{2}$.

Following the universal character of PLE for all living systems and the starting hypothesis of this work that living agents can be considered as thermal engines, we hope that this approach is meaningful to most if not all living systems and systems used or driven by living agents. We can cite, among others, economic systems, linguistic systems, informational and social networks and so on, which are systems driven by the effort of human beings to achieve something. For example, in an economic system, income is the achievement $X$ people are trying to obtain with the effort. People will in general try to achieve more income with less effort, producing the Pareto law in the distribution of income [15]. In the linguistic systems, a similar thing happens. People try not to be garrulous (statistically, of course) and to express as much information as possible using as less words as possible. This behavior tends to increase the frequency of individual words in a text of a given number of words, as if increase of frequency of words was the thing people are looking for with effort. This produces the Zipf law of frequency-rank distribution [2][3].

If we think about city populations in a society, people tend to move to large cities for many reasons, as if the increasing size of city was something they were trying to obtain with of course as little effort as possible. This produces the Zipf and Pareto laws of size distribution [4][16]. The last example we think worth mentioning concerns the informational networks where the Zipf and Pareto laws are a consequence of the rule of preferential attachment [17]. This rule is a typical behavior of least effort or of maximum efficiency in the quest for information. Its derivation from PLE will be discussed in a coming work [18].

We would like to remark that the Zipf and Pareto laws are often related to the Pareto 80/20 rule[4], or same rule with a slight different proportion (90/10 for instance), in the sense that this rule can be calculated from the Pareto law. The Pareto rule is easy to observe, it is an interesting illustration of the ubiquity of the Zipf and the Pareto laws, you can hardly ignore it even in daily life. For example, if you have a large number of phone numbers in your calling list, you can easily check that 80% (or even more) of your call are made to only 20% (or even less) of the phone numbers in that list. The same phenomenon occurs to your mailing list. If you have many friends, you can easily check that a small percentage (perhaps 10%) of them take up most (perhaps 90%) of your time spent with friends. Similar phenomena happen with the books you keep in your library, and your shoes and clothes. All these phenomena illustrate the ubiquity of the Zipf and the Pareto laws as well as of the underlying PLE.

### 4) Uniqueness of the efficiency functional

An important step of our application of PLE to derive power laws is to discover the efficiency functional Eq.(5) from the consideration of the independent probabilities of subsytems Eq.(3) and the nonadditive property Eq.(4) of efficiency. The efficiency functional was proposed in [14] in an intuitive way. It is very important to prove that Eq.(3) and Eq.(4) uniquely lead to Eq.(5). This uniqueness is a guarantee of the reliability of the efficiency functional as a measure of performance and a maximizable fundamental quantity to yield power laws. Below we provide a rigorous proof of the uniqueness of the efficiency functional of Eq.(5).

The uniqueness of Eq.(5), or $(1 + a\eta_i) = p_i^b$, for given Eq.(3) and Eq.(4), can be proved by using the following lemma.

**Lemma:**

---

[4] https://en.wikipedia.org/wiki/Pareto_principle

If $f$ is a positive continuous function on the σ-algebra $B(R^+)$ satisfying the condition

$$\forall x \in \mathbb{R}^+, \forall y \in \mathbb{R}^+, f(xy) = f(x)f(y),$$

then there is $a > 0$ such that $\forall x \in \mathbb{R}^+, f(x) = x^a$.

**Proof**

Let $\forall x \in \mathbb{R}^+, g(x) = \ln(f(e^x))$, $g$ being continuous on $\mathbb{R}^+$.

We have for every $x \in \mathbb{R}^+$ and $y \in \mathbb{R}^+$, the following relationship:

$g(x+y) = \ln(f(e^{x+y})) = \ln(f(e^x e^y)) = \ln(f(e^x)f(e^y))$ following the hypothesis in the lemma. We get $g(x+y) = g(x) + g(y)$.

Consequently, $\forall n \in \mathbb{N}^*, g(n) = ng(1), g(1) = g\left(\frac{1}{n} + \cdots + \frac{1}{n}\right) = ng\left(\frac{1}{n}\right)$, et $g\left(\frac{1}{n}\right) = \frac{g(1)}{n}$.

Thus, $\forall p \in \mathbb{N}^*, \forall q \in \mathbb{N}^*, g\left(\frac{p}{q}\right) = g\left(\frac{1}{q} + \cdots + \frac{1}{q}\right) = pg\left(\frac{1}{q}\right) = p\frac{g(1)}{q}$.

Because for $\forall x \in \mathbb{R}^+$, there is a sequence $(p_n, q_n)$ in $\mathbb{N}^* \times \mathbb{N}^*$, such that $\lim_{n \to \infty} \frac{p_n}{q_n} = x$. On the other hand, as $g$ is continuous in $x$, we can write

$$\lim_{n \to \infty} g\left(\frac{p_n}{q_n}\right) = g(x) = \lim_{n \to \infty} p_n \left(\frac{g(1)}{q_n}\right) = xg(1)$$

Now, $\forall x \in \mathbb{R}^{*+}$, let $X = \ln x$, this gives

$$f(x) = f(e^X) = e^{g(X)} = e^{Xg(1)} = (e^X)^{g(1)} = x^a,$$

with $a = g(1)$. The uniqueness of the relation $(1 + a\eta_i) = p_i^b$ and of Eq.(5) is thus proven.

## 5) Lesche stability of the efficiency functional

Now let us look at a mathematical property of the efficiency functional Eq.(5) concerning its robustness against small variation in the probability distributions. The idea is the following: if a continuous function represents a physical quantity which is a function of some variables such as time, position, energy, probability distribution etc., the definition of the function is robust (or Lesche stable) if and only if the function undergoes smooth variation when the variables smoothly change due to the evolution of the system [19]. This robustness guarantees that the defined function is a valid and meaningful physical quantity. This condition has been used as a criterion to verify the robustness of several generalizations of entropy, information and statistical means. The conclusion is that not all continuous and concave functionals are Lesche stable or robust against small variation of probability [19][20]. In what follows, we provide a proof of the robustness (Lesche stability) of the efficiency functional Eq.(5) to guarantee the efficiency functional as a valid and meaningful measure of performance.

**Definition.** A Function $C$ defined on $A = \bigcup_{n \in \mathbb{N}^*}\{p \in [0,1]^n, \sum_{i=1}^n p_i = 1\}$ is stable (robust) if $C$ has the following property

$\forall \varepsilon > 0, \exists \delta > 0, \forall N \in \mathbb{N}^*, \forall p, p' \in A \cap \mathbb{R}^N$, such that $\|p - p'\|_1 < \delta \Rightarrow \left|\frac{C(p) - C(p')}{C_{N,max}}\right| < \varepsilon$

where $C_{N,max}$ is defined by $C_{N,max} = \max\{|C(p)|, p \in A \cap \mathbb{R}^N\}$.

**Proposition.**

For all $a \in \,]0,1[$, the efficiency $p \mapsto E(p)$ given by Eq.(5) is stable.

To establish the proof of this proposition, we need the following lemmas.

**Lemma 1.**
Let $a \in \,]0,1[$, $N \in \mathbb{N}^*$, $p = (p_i)$, $p' = (p_i')  \in [0,1]^N$. We have

$$\sum_{i=1}^N \left|p_i^{1-a} - (p_i')^{1-a}\right| \leq N^a (\|p - p'\|_1)^{1-a}$$

**Proof**:

Suppose $p \neq p'$; let

$$I_1 = \{i \in \{1, \cdots, N\}, p_i' < p_i\} \text{ and } I_2 = \{i \in \{1, \cdots, N\}, p_i < p_i'\}.$$

We have $I_1 \neq \emptyset$, $I_2 \neq \emptyset$, and for all $i \in I_1$,

$$1 = \left(\frac{p_i'}{p_i}\right) + \left(1 - \frac{p_i'}{p_i}\right) \leq \left(\frac{p_i'}{p_i}\right)^{1-a} + \left(1 - \frac{p_i'}{p_i}\right)^{1-a} = \left(\frac{p_i'}{p_i}\right)^{1-a} + \left|1 - \frac{p_i'}{p_i}\right|^{1-a}.$$

We can deduce $(p_i)^{1-a} - (p_i')^{1-a} \leq |p_i - p_i'|^{a1-a}$. By symmetry, we also get $(p_i')^{1-a} - (p_i)^{1-a} \leq |p_i - p_i'|^{1-a}$ for all $i \in I_2$.

It can be concluded that for $i \in \{1, \cdots, N\}$, $\left|(p_i)^{1-a} - (p_i')^{1-a}\right| \leq |p_i - p_i'|^{1-a}$.

Now, we write $\sum_{i=1}^N |(p_i)^{1-a} - (p_i')^{1-a}| \leq \sum_{i=1}^N |p_i - p_i'|^{1-a} = \sum_{i=1}^N 1 \times |p_i - p_i'|^{1-a}$. Considering Holder's Inequality $\left(\text{because } 1 = a + (1-a) = \frac{1}{1/a} + \frac{1}{1/(1-a)}\right)$, we obtain

$$\sum_{i=1}^N 1 \times |p_i - p_i'|^{1-a} \leq \left(\sum_{i=1}^N 1^{\frac{1}{a}}\right)^a \left(\sum_{i=1}^N \left(|p_i - p_i'|^{1-a}\right)^{\frac{1}{1-a}}\right)^{1-a}$$

$$\leq N^a \left(\sum_{i=1}^N |p_i - p_i'|\right)^{1-a}$$

$$= N^a (\|p - p'\|_1)^{1-a}$$

**Lemma 2.**

For all $a \in\ ]0,1[$, the function $x \mapsto \frac{x^a}{|1-x^a|}$ is bounded on $[2,+\infty[$ by $M > 0$.

Indeed if $a \in\ ]0,1[$, the function $x \mapsto \frac{x^a}{|1-x^a|}$ is continuous on $[2,+\infty[$ and has $\lim\limits_{x \to +\infty} \frac{x^a}{|1-x^a|} = 1$.

**Proof of the proposition.**

Using the Lagrange multiplier method, we prove that:

$$E_{N,max} = E\big((1/N),\cdots,(1/N)\big) = \frac{|1-N^a|}{a}.$$

Let $\varepsilon > 0$, there exists $\delta = \left(\frac{\varepsilon}{M}\right)^{1/(1-a)} > 0$, such that for all $N \in \mathbf{N}^*\backslash\{1\}$, $p, p' \in A \cap \mathbb{R}^N$, $\|p - p'\|_1 < \delta$ implies

$$\begin{aligned}
\left|\frac{E(p) - E(p')}{E_{N,max}}\right| &= \frac{\left|\sum_{i=1}^N \left(p_i^{1-a} - (p_i')^{1-a}\right)\right|}{|1 - N^a|} \\
&\leq \frac{\sum_{i=1}^N \left|p_i^{1-a} - (p_i')^{1-a}\right|}{|1 - N^a|} \\
&\leq \frac{N^a\left(\sum_{i=1}^N |p_i - p_i'|\right)^{1-a}}{|1 - N^a|} \\
&= \frac{N^a(\|p - p'\|_1)^{1-a}}{|1 - N^a|} \\
&< M\delta^{1-a} = \varepsilon
\end{aligned}$$

The above proposition is thus proved, meaning that the definition Eq.(5) of the efficiency is robust and Lesche stable against small variations or perturbations of the probability. As mentioned above, this is a basic guarantee for this functional to be a valid physical quantity and a reliable measure of efficiency for the application of PLE and for the metric of performance.

## 6) Negative ZP efficiency?

In the work [14], the efficiency functional Eq.(5) was supposed to be positive in the case of discrete distribution. Negative efficiency in thermodynamics would mean an engine absorbing heat but doing negative work (consuming work), which is not meaningful from the viewpoint of physics. This positivity is guaranteed by Eq.(5) only if $X$ is a discrete variable. Now if $X$ becomes continuous with a continuous probability distribution $\rho(x)$, Eq.(5) becomes

$$\eta = \int_{x_{min}}^{x_{max}} \frac{\rho(x)^{1-a} - 1}{a} dx \tag{6}$$

Negative efficiency turns out to be possible depending on the distribution (see below). This may not be possible for a thermal engine, but may make sense in many domains outside physics if the achievement is negative. In economics for example, if a person or a company is indebted during a given period, negative efficiency is useful to describe the reality. The analysis below allowing negative efficiency is motivated by this consideration.

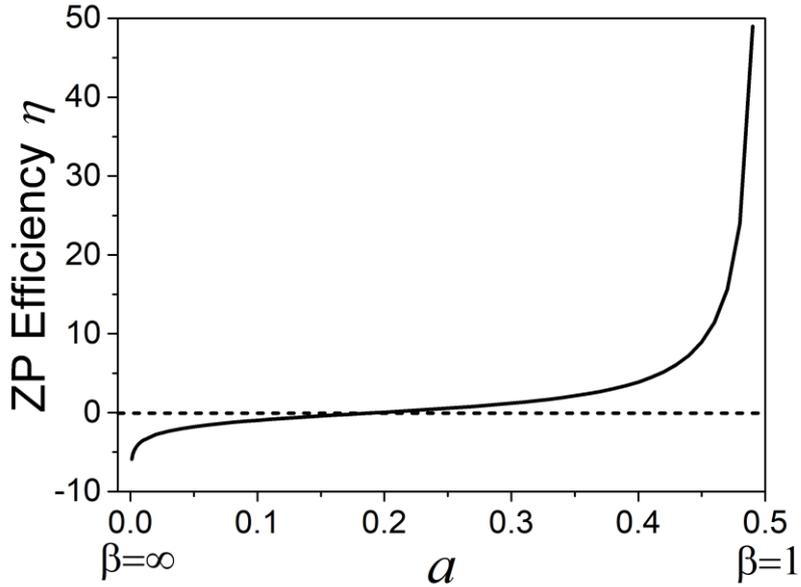

Figure 2. Evolution of the efficiency of Eq.(6) calculated for Pareto's distribution in the interval $0 < a < 0.5$ or $\infty < \beta < 1$. $\eta$ diverges for $a=0.5$ ($\beta = 1$) and $a=0$ ($\beta = \infty$). $\eta = 0$ for $a = 0.1945$ or $\beta = 4.14$, and $\eta < 0$ in the range $0 < a < 0.1945$ or $4.14 < \beta < \infty$.

As shown in [14], the PLE implies the maximization of the efficiency Eq.(6), which leads to Pareto's law

$$P(X > x) = \left(\frac{x_m}{x}\right)^\beta \qquad (7)$$

where $P(X > x)$ is the probability of finding a person with income larger than a value $x$, $x_m$ is the smallest income and $\beta = \frac{1}{a} - 1$ is a positive constant characterizing the distribution. We have $0 < a < 1$ and $0 < \beta < \infty$. With this distribution, the efficiency Eq.(6) can be calculated with the Pareto probability density distribution (PDF) $\rho(x) = \frac{\beta}{x^{\beta+1}}$ where we suppose $x_m = 1$ and $x_{max} = \infty$ for simplicity. The result reads $\eta = \frac{1}{a}\left[\left(\frac{a}{1-a}\right)^a \frac{1-a}{1-2a} - 1\right] = (\beta + 1)[\beta^{-\frac{1}{\beta+1}} \cdot$ $\frac{\beta}{\beta-1} - 1$] which is referred to as Zipf and Pareto (ZP) efficiency and plotted in Figure 2 for $0 < a < 0.5$. We have excluded $0.5 < a < 1$ due to the divergent mean of $x$ [14].

The ZP efficiency $\eta$ increases with increasing $a$ or decreasing $\beta$, and diverges for $a=0.5$ ($\beta = 1$) and $a=0$ ($\beta = \infty$). $\eta = 0$ for $a = 0.1945$ or $\beta = 4.14$. $\eta < 0$ for $0 < a < 0.1945$ or $4.14 < \beta < \infty$. The increase of ZP efficiency with decreasing $\beta$ (or increasing $a$) is quite natural because, considering the Pareto PDF $\rho(x) = \frac{\beta}{x^{\beta+1}}$, decreasing $\beta$ implies increasing probability of large achievement $x$ and increasing global achievement represented by the average achievement $\bar{x} = \frac{1}{1-1/\beta}$. The divergence of the ZP efficiency and of $\bar{x}$ at $\beta = 1$ is of course due to the choice of the infinite upper limit of $x$ ($x_{max} = \infty$). If $x_{max}$ is finite, the ZP efficiency should be finite, but can be very large at $\beta = 1$.

## 7) Efficiency and inequality

The relationship between efficiency and inequality is a hot topic in economics [21][22][23]. From a classical point of view, economic efficiency or growth is positively correlated with inequality. This belief comes from the idea that inequality is a source of incentive. Inequality is a motivator and can incite effort to raise income, invest, work hard, and innovate. Without opportunity, firms and individuals will reduce effort resulting in lower economic growth. The economic disaster in some former communist countries which carried out egalitarian experiments is an example of the failure of redistributive policies. The same reasoning makes sense in educational systems as well. Inequality in educational levels implies, at least statistically, opportunity to receive higher education, to attain more knowledge, higher social class, higher income, honor and accolades and so on. It is imaginable that if everybody in a population has the same educational level, the effort to achieve more education, more knowledge and skills will be reduced or even non-existent, leading to low efficiency marked by low production of knowledge and innovation. This incentive mechanism has its place in every system of living agents where effort is needed to obtain higher achievement.

However, a question naturally arises as to what extent the incentive of inequality works as a driving force. It is not convincing to claim that an economic system is extremely efficient when only one person has all the wealth and the rest of the population has nothing. One cannot say that an educational system has maximal efficiency when only one person is very educated and the rest of the population has no education at all (corresponding to a Gini coefficient $G=1$ [24]). The reality is much more complex than the role of inequality being a pure incentive . We can get an idea about this complexity from a figure of the evolution of Gini coefficients after World War II available on a website [24]. In this figure, we see many countries (USA, China, India, UK etc.) who experienced increasing inequality accompanied by long-term growth. On the other hand, we see obvious counterexamples with countries such as France, Germany, Japan etc. showing decreasing inequality coefficient during the period of long-term growth. The existence of different theories of economic growth is a reflection of the diversity of viewpoints in this domain [23][24].

Actually, economic growth depends on so many non-economic or external factors like political and social stability, technological level and progress, human capital and so on, that sometimes the intrinsic factors such as incentive may be hidden, making it difficult to figure out how and to what degree the economic performance depends on the intrinsic factors. For this reason, it is of interest to study the behavior of the ZP efficiency in relation with the inequality measured by the Gini coefficient $G$. For a system obeying the Pareto law, the Gini coefficient reads $G = \frac{1}{2\beta - 1}$. In Figure 3 we see a positive correlation between the ZP efficiency $\eta$ and the Gini coefficient $G$. $\eta$ becomes very large as $G \to 1$, which is the case of maximal inequality. $\eta$ decreases with decreasing $G$ and becomes negative when $G$ is smaller than a given value corresponding to the Pareto index $\beta = 4.14$.

Usually, the inequality $G$ should be considered as an intrinsic property of the system of interest [21]-[24]. So from the positive correlation between the ZP efficiency and the inequality, we can say that the ZP efficiency should be regarded as an intrinsic property of the system. This implies that its behavior in Figure 2 and Figure 3 is an intrinsic behavior of the system of interest without the influences of the external factors. We would say that this attribute of the ZP efficiency is quite natural because this efficiency is derived from the pure mathematical property of nonadditivity Eq.(4), which is a system-independent property of any ratio between output/achievement and input/effort [27]. This insight into the ZP efficiency is important for it to be used as a general metric of performance. To our eyes, this intrinsic and system-independent attribute of the ZP efficiency is a guarantee of its reliability as a general measure of performance. Despite its pure mathematical character, it should be able to reflect, indirectly, the complexity of the interaction and interplay among the intrinsic and extrinsic factors in a complex system when it is calculated from the knowledge of the probability distribution laws of the achievement $X$. These laws, derived from empirical investigation or numerical simulation, carry important information about the functioning of the system. Having said this, we agree that the question whether the ZP efficiency can be used as a practical metric of performance should be investigated further, with the help of empirical data as well as numerical simulation.

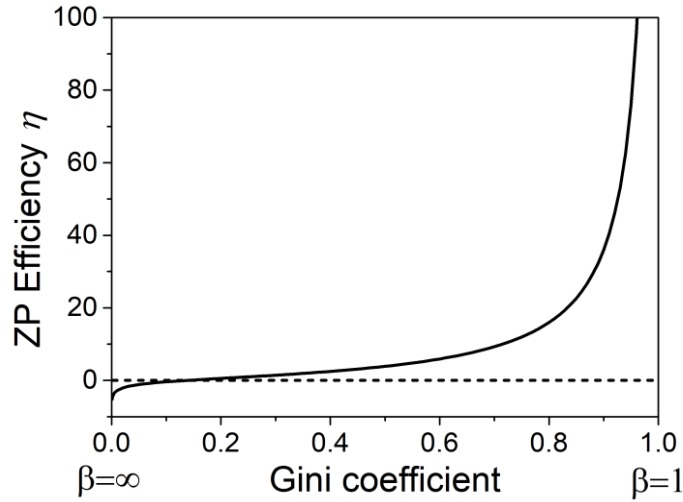

Figure 3. Variation of the ZP efficiency $\eta$ as a function of the Gini coefficient $G$ indicating inequality ($G$=0 and $G$=1 implies absolute equality and maximum inequality, respectively). $\eta$ becomes negative in the interval $0 < G < 0.14$ corresponding to $0 < a < 0.1945$ or $4.14 < \beta < \infty$, and very large when $G \to 1$ (maximum inequality).

## 8) Relationship with some entropies

In order to get more insight into the efficiency Eq.(5) or (6), it is of interest to compare it to the informational entropy of the same probability distribution laws as used in the calculation of the efficiency. This informational entropy is another important quantity characterizing the state of the systems concerning the uncertainty or disorder in its distribution of achievement. In what follows, we present the results using a generalized uncertainty measure called varentropy [25][26]. It is worth mentioning that many widely used entropy formulas, including the Shannon entropy, are not suitable for continuous probability distribution for several mathematical reasons widely discussed in [29][31][32][33]. One of these reasons is the problem of negative entropy (see [29][34] and references therein). There are some solutions using for example Kullback-Leisler relative entropy proposed in [29][31][32][33]. But the measure of the absolute uncertainty of continuous probability distribution, in the same way as for discrete probability, is still a topic of discussion. We propose a pathway here using varentropy to measure probabilistic uncertainty. Varentropy is defined by mimicking the variational form of the entropy of the second law of thermodynamics [25]. This is a definition from scratch without any prerequisite or postulate on the property of entropy. The motivations was to seek an uncertainty measure that is maximal or maximizable for any probability distribution. Although $S_{BS}$ has been widely used as a universal uncertainty measure for any probability distribution, it is maximizable only for exponential and uniform distribution [25][30]. In other words, the outcome of its maximization by the calculus of variation is only the exponential distribution or the uniform distribution. This fact was especially noticed during the development of several generalized statistics based on the principle of maximum entropy [8]-[12]. During this development, many informational entropy formulas were proposed, each one being maximized by a specific probability distribution. A question was thus raised [25][30] about whether it was possible to define a general informational entropy (or measure of uncertainty) from the point of view of physics that can be maximized to yield any probability distribution. In other words, it should be just $S_{BS}$ maximized by exponential distribution, and other known entropies which are maximized by their own probability distributions, respectively [25][26][28][30].

Such a generalized measure does exist [25]. It is given by $\delta S_V = A(\delta \bar{x} - \overline{\delta x}) = A \int x \, \delta \rho(x) dx$ for a continuous distribution $\rho(x)$, where $S_V$ is the uncertainty measure named varentropy, $\bar{x}$ is the mean of the random variable $x$, $\overline{\delta x}$ is the mean of a small variation of $x$, and $A$ is a constant determining the dimension and the sign of the uncertainty measure. When the variable $x$ is energy $E$ of a thermodynamic system, we have $\delta S_V = A(\delta \bar{E} - \overline{\delta E}) = \frac{1}{T}(\delta U - \delta W)$, where $A = \frac{1}{T}$ the inverse temperature, $\bar{E} = U$ the internal energy, $\delta W = \overline{\delta E}$ is the external work and $S_V$ is just the entropy of the second law of thermodynamics. The maximization of $S_V$ is then a natural consequence of the second law of thermodynamics. This entropy turns out to be equivalent to a generalized entropy defined in [30] on the basis of the principle of maximum entropy. It has been shown that, for discrete probability, varentropy is just the Shannon entropy maximized by exponential distribution, the Tsallis (or Renyi) entropy maximized by q-exponential distribution, the power law entropy [25] maximized by the Zipf and Pareto laws, and many other entropy formulas maximized by their own probability distributions, respectively. The reader is referred to [25][26][28][30] for more details about these entropies and their probability distributions.

It is worth mentioning that, for several probability distributions different from the exponential distribution, varentropy is always larger than the Shannon entropy. An example is in Figure 4 below where the varentropy maximized by Pareto law is shown to be equal to or larger than the Shannon like informational entropy (maximized by the exponential distribution) calculated for the same power law. Another example was given for a double power law in [28]. Actually, varentropy is expected to be larger than (or equal to) the Shannon entropy for any probability distribution different from the exponential one, because varentropy is always maximized for whatever distribution but Shannon entropy is only maximized by the exponential distribution. In this sense, varentropy is a more optimal metric of probabilistic uncertainty and disorder.

In this work, we focus on another advantage of varentropy which allows avoiding negative entropy for continuous distributions. It is well known that Shannon entropy (and many others including the Tsallis and Renyi entropies) is by definition always positive for discrete distributions, but can become negative when applied to continuous distributions. Some examples of negative entropy given by the Shannon formula for the continuous exponential and power law distributions can be seen in [34] and references therein.

Now let us look at the varentropy for continuous probability distributions. For the continuous exponential distribution $\rho(x) = \frac{1}{Z} e^{-x}$ where $x$ is a dimensionless variable and Z the normalization constant, we have, by definition, $\delta S_V = -A \int \ln(Z\rho) \, \delta\rho dx = -A \int (\ln Z + \ln \rho) \, \delta\rho dx = -A\delta \int \rho \ln \rho \, dx = \delta[-A \int \rho(\ln \rho - \ln m) \, dx]$ where $m(x)$ is some function satisfying $\int \rho \ln m(x) \, dx = C$, where $C$ is a constant with $\delta C = 0$. Let A=1, we obtain the continuous Boltzmann-Shannon like entropy denoted by $S_V^{BS}$ in keeping the index $V$ to indicate its origin from varentropy:

$$S_V^{BS} = -\int \rho \ln \frac{\rho}{m} dx \tag{8}$$

where $m$ is the invariant measure similar to that introduced by Jaynes [29].

The formula Eq.(8) can take another equivalent form if we keep the term lnZ in the above calculation, i.e. $\delta S_V^{BS} = -A \int (\ln Z + \ln \rho) \, \delta\rho dx = -A\delta \int \rho(\ln Z + \ln \rho) \, dx = \delta[-A \int \rho(\ln(Z\rho) - \ln m) \, dx]$ which implies $S_{BS} = -\int \rho \ln \frac{Z\rho}{m} dx$ (A=1).

Eq.(8) can be used to measure the uncertainty in the Pareto distribution $\rho(x) = \frac{\beta}{x^{\beta+1}}$. It is easy to calculate from Eq.(8): $S_V^{BS} = 1 + \frac{1}{\beta} + C$. Let $C = -1$, we get

$$S_V^{BS} = \frac{1}{\beta} = \frac{a}{1-a} \geq 0 \tag{9}$$

which is plotted in Figure 4.

It is noteworthy that, as well known, $S_V^{BS}$ can be maximized to obtain uniform and exponential distributions. Hence, it is not at maximum for power law distributions. For a discrete power law $p_i = \frac{1}{Z} x_i^{-\frac{1}{b}}$, the maximized entropy is $S_V = \frac{\sum_i p_i^{1-b}-1}{1-b}$ which is positive and zero for non-probabilistic case [25]. This entropy looks like Tsallis entropy $S_T = \frac{\sum_i p_i^q - 1}{1-q}$ but is essentially different. $S_T$ is a generalized $q$-logarithmic entropy (while $S_V$ is not) and tends to Shannon entropy when $q \to 1$, but $S_V$ tends to zero for $b \to 0$ and diverges for $b \to 1$ [25].

Now let us calculate the varentropy $S_V^P$ for continuous power law distribution $\rho(x) = \frac{1}{Z} x^{-\frac{1}{b}}$ directly from its definition: $\delta S_V^P = A \int (Z\rho)^{-b} \delta\rho dx = A \int \delta \left( \frac{Z^{-b}}{1-b} \rho^{1-b} \right) dx = \delta \{ \frac{A}{1-b} \int (Z^{-b} \rho^{1-b} - m) dx \} \}$. Let $A = 1$, this continuous varentropy reads

$$S_V^P = \int \rho \frac{(Z\rho)^{-b} - m}{1-b} dx \tag{10}$$

where the function $m$ is such that $\int \rho m(x) dx = C$ is a constant of the variation, i.e., $\delta C = 0$.

For Pareto PDF $\rho(x) = \frac{\beta}{x^{\beta+1}}$, let $Z = \frac{1}{\beta}, b = \frac{1}{\beta+1}$ in Eq.(12), we obtain

$$S_V^P = \frac{\int_1^{+\infty} \beta^{\frac{1}{\beta+1}} \cdot \rho^{\frac{\beta}{\beta+1}} dx - C}{\frac{\beta}{\beta+1}} = \frac{\beta+1}{\beta} \left( \frac{\beta}{\beta-1} - C \right). \tag{11}$$

Let $C = 1$, we obtain a positive varentropy

$$S_V^P = \frac{\beta+1}{\beta(\beta-1)} \geq 0 \tag{12}$$

This continuous varentropy of the Pareto distribution is plotted in Figure 4 in comparison with the ZP efficiency $\eta$ and the continuous $S_V^{BS} = \frac{1}{\beta} = \frac{a}{1-a}$ calculated from the Pareto distribution. We observe that the efficiency is positively correlated to entropy as a measure of dynamical uncertainty or disorder. We also noticed that $S_V^P$ is always larger than $S_V^{BS}$, which is expected because $S_V^P$ is the varentropy maximized by the power law itself, so it is more optimal in measuring the power law uncertainty than $S_V^{BS}$ which is only optimal or maximized for exponential distributions.

By the way, this result also reveals a positive correlation between inequality and dynamic disorder, implying that larger inequality occurs with larger uncertainty or disorder in the probability distribution of achievement. The full understanding of this correlation needs further investigation.

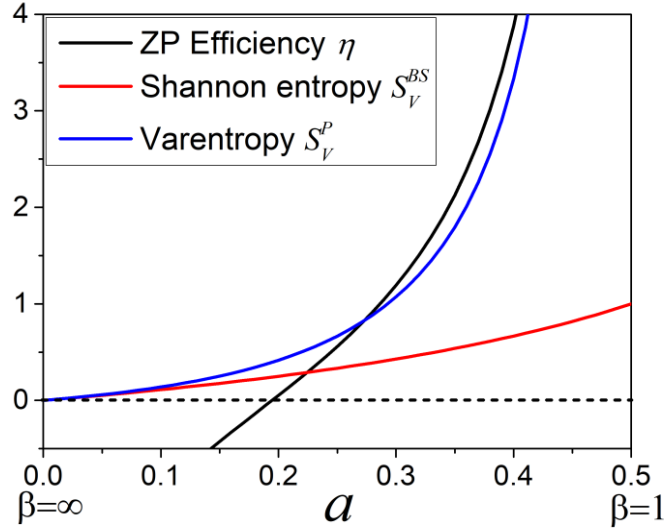

Figure 4. Evolution of ZP efficiency $\eta = (\beta + 1)[\beta^{-\frac{1}{\beta+1}} \cdot \frac{\beta}{\beta-1} - 1]$ compared to the continuous $S_V^{BS} = \frac{1}{\beta} = \frac{a}{1-a}$ and $S_V^P = \frac{\beta+1}{\beta(\beta-1)}$. Both entropies increase from zero with increasing $a$ or decreasing $\beta$ in the intervals $0 < a < 0.5$ or $1 < \beta < \infty$. So they have positive correlation with the ZP efficiency. However, at $a \to 0.5$ or $\beta \to 1$, $S_V^{BS}$ is finite and $S_V^P$ tends to infinity. So $S_V^P$ is closer to $\eta$ than $S_V^{BS}$. Notice that $S_V^P \geq S_V^{BS}$ as expected.

## 9) Conclusion

In this paper, an efficiency functional used in the derivation of the Zipf and Pareto laws from the principle of least effort has been investigated in order to get further insight into its mathematical properties and its potential as a general measure of performance. We have proved that this efficiency functional is the correct one and that it is unique under the condition of the nonadditivity of efficiency of thermodynamics. We have also provided a proof of its robustness against small variations of probability. These two mathematical properties of the efficiency functional guarantee its validity and reliability as a measure of performance and an origin of the power laws of complex living systems.

The behavior of this efficiency has been studied for the Pareto distribution and compared to the Gini coefficient (inequality) and probabilistic uncertainty measured by a generalized informational entropy called varentropy. Varentropy is an optimal measure of probabilistic uncertainty and allows us to avoid negative value of uncertainty measure. We have shown that the value of the efficiency increases with decreasing index *β* of the Pareto law, implying that the efficiency is positively correlated with the statistical mean of the achievement. This result confirms the usefulness of the efficiency functional as a metric of performance.

It has also been shown that the Zipf-Pareto efficiency is positively correlated with the Gini coefficient and probabilistic uncertainty (entropy). To put it differently, the Zipf-Pareto efficiency increases with increasing inequality and probabilistic uncertainty in the system of interest. This result also implies that inequality has a positive correlation with the probabilistic uncertainty (disorder). In other words, the inequality increases with increasing disorder in the system.

As mentioned above, the objective of this work is, after the previous work using the efficiency functional to derive the Zipf and the Pareto laws from PLE, to investigate this functional from a mathematical point of view in order to better understand it as a metric of performance of the systems showing the Zipf and the Pareto laws. The present work has confirmed the usefulness of this efficiency for these systems. Nevertheless, the question whether or not this efficiency can be a useful measure of performance for other systems showing more complicated distribution laws is still open. This question is related to another open question which is whether it is possible to derive more complicated distribution laws from the principle of least effort. Further work is in progress in this direction.

The principle of least effort really seems a universal fundamental rule for the world of the living according to the analysis of Ferrero who thought that the origin of this rule was the laws of inertia in physical and chemical motions [1]. It is true that a living agent should always try to economize its energy, to achieve more by doing less, or to achieve the necessary things by doing the minimum (see principle of good enough [35] for example). This is even a necessary condition for its survival and for the existence of its species. Our question is how far we can go with this universal principle in the study of the behaviors, the empirical laws, the past and the future of the living world.

## Data availability

The data that support the findings of this study are available from the corresponding author upon reasonable request.